\documentclass[12pt, reqno]{amsart}

\usepackage{iftex}
\ifPDFTeX
  \usepackage[T1]{fontenc}
  \usepackage[utf8]{inputenc}
  \IfFileExists{lmodern.sty}{\usepackage{lmodern}}{}
\else
  \usepackage[no-math]{fontspec}
\fi
\usepackage[l2tabu,orthodox]{nag}
\usepackage{amsmath, amssymb, amsthm}
\usepackage{mathtools}
\usepackage[all,error]{onlyamsmath}
\usepackage[strict]{csquotes}
\usepackage{mathrsfs}
\usepackage[dvipsnames]{xcolor}
\usepackage{cite}
\usepackage[colorlinks, linkcolor=red, citecolor=blue, urlcolor=blue, hypertexnames=false]{hyperref}

\theoremstyle{plain}
\newtheorem{thm}{Theorem}[section]
\newtheorem{lem}[thm]{Lemma}
\newtheorem{prop}[thm]{Proposition}
\newtheorem{cor}[thm]{Corollary}
\newtheorem{conj}[thm]{Conjecture}
\theoremstyle{definition}
\newtheorem{defi}[thm]{Definition}
\newtheorem{prob}[thm]{Open problem}
\theoremstyle{remark}
\newtheorem{rem}[thm]{Remark}

\numberwithin{equation}{section}

\newcommand{\R}{\mathbb{R}}

\newcommand{\mS}{\ensuremath{\mathbb{S}}}

\newcommand{\jb}[1]{\langle #1\rangle}
\newcommand{\dd}{\mathrm{d}}

\author[R. Castillo]{Ricardo Castillo}
\address[R. Castillo]{Departamento de Matem\'atica, Universidad del B\'io-B\'io (UBB), Concepci\'on, Chile.}
\email{rcastillo@ubiobio.cl}
\urladdr{https://orcid.org/0000-0002-0034-2006}
\author[M. Majdoub]{Mohamed Majdoub}
\address[M. Majdoub]{Department of Mathematics, College of Science, Imam Abdulrahman Bin Faisal University, P. O. Box 1982, Dammam, Saudi Arabia.}
\address[M. Majdoub]{Basic and Applied Scientific Research Center, Imam Abdulrahman Bin Faisal University, P.O. Box 1982, 31441, Dammam, Saudi Arabia.}
\email{mmajdoub@iau.edu.sa}
\email{med.majdoub@gmail.com}
\email{mohamed.majdoub@fst.rnu.tn}
\urladdr{https://orcid.org/0000-0001-6038-1069}
\title[Grow-up for forced viscous Hamilton--Jacobi equations]{Growth versus boundedness for viscous Hamilton--Jacobi equations with a forcing term: a Fujita-type dichotomy}

\date{\today}

\begin{document}



\begin{abstract}
We study the Cauchy problem for $u_t=\Delta u+|\nabla u|^q+\mathtt a(t)\mathbf w(x)$ in $(0,\infty)\times\R^N$, where $q>1$ and the forcing profile $\mathbf w\in L^1(\R^N)$ satisfies $\int\mathbf w>0$. In contrast with the semilinear problem $u_t=\Delta u+|u|^p+\mathbf w(x)$ studied by Bandle, Levine and Zhang, solutions with regular data never blow up in finite time. The natural analogue of the Fujita problem is therefore to decide whether solutions remain bounded or grow up as $t\to\infty$. A scaling argument and an elliptic Liouville theorem single out the critical exponent $q_c=N/(N-1)$, with $q_c=+\infty$ when $N=1$. For time-independent forcing, we prove that every solution grows at least linearly in time on a fixed ball when $1<q<q_c$, and also at the critical exponent $q=q_c$ when $N\ge2$. In the plane, every solution tends to $+\infty$ locally uniformly for every $q>1$, at least at the rate $(4\pi)^{-1}\big(\int\mathbf w\big)\log t$. In the quadratic case $q=2$, the Cole--Hopf transformation shows that $\sup u(t)$ grows at least like $-\lambda_0t$, where $\lambda_0$ is the bottom of the spectrum of $-\Delta-\mathbf w$. Known results on Schr\"odinger semigroups show that this rate is exact. In particular, the growth is linear when $N\le2$, even for nontrivial forcing with zero mass. Conversely, if $N\ge3$ and $q>q_c$, solutions with small nonnegative data and forcing remain bounded for all times. For time-dependent profiles, we set $A(t)=\int_0^t\mathtt a(s)\,\dd s$. We show that solutions grow like $A(t)$ whenever $A$ grows superlinearly, and that they grow at least linearly on a fixed ball when $1<q\le q_c$ and $\liminf_{t\to\infty}A(t)/t>0$. We also give a partial result for decaying profiles. Several conjectures and open problems are formulated, and proved statements are kept separate from heuristics.
\end{abstract}

\maketitle

\section{Introduction}\label{sec:intro}

We consider the Cauchy problem
\begin{equation}\label{eq:main}
\begin{cases}
u_t=\Delta u+|\nabla u|^q+\mathtt a(t)\,\mathbf w(x), & (t,x)\in(0,\infty)\times\R^N,\\
u(0,x)=u_0(x), & x\in\R^N,
\end{cases}
\end{equation}
where $q>1$, $\mathtt a$ is a continuous time profile (in most of the paper $\mathtt a\equiv1$), and $\mathbf w$ is a forcing profile with positive mass, that is,
\begin{equation}\label{eq:mass}
\mathbf w\in L^1(\R^N)\cap W^{1,\infty}(\R^N),\qquad M_\infty:=\int_{\R^N}\mathbf w(x)\,\dd x>0 .
\end{equation}
The sign of $\mathbf w$ is {not} prescribed. Only its total mass is assumed to be positive.

Equation \eqref{eq:main} is a viscous Hamilton--Jacobi equation with a source term. For $q=2$, it is the deterministic version of the Kardar--Parisi--Zhang equation, which describes the growth of an interface by deposition \cite{KPZ86}. In this model, $u$ is the height of the interface, the term $\Delta u$ accounts for surface relaxation, the term $|\nabla u|^2$ for lateral growth, and $\mathtt a(t)\mathbf w(x)$ is the local deposition rate. General exponents $q$ appear in deterministic models of surface growth \cite{KS88}. From this point of view, the question studied below is natural. {Does a localised deposition rate with positive mass make the interface grow without bound, or can the interface settle to a bounded profile?}

Our motivation comes from the theory of critical exponents for forced parabolic problems, which goes back to the work of Fujita on $u_t=\Delta u+u^p$ \cite{fujita1966}. For the problem $u_t=\Delta u+|u|^p+\mathbf w(x)$ with $\int\mathbf w>0$, Zhang \cite{Zha98} and Bandle, Levine and Zhang \cite{BLZ} proved the following dichotomy. If $N\le2$, or if $N\ge3$ and $1<p\le N/(N-2)$, then no global solution exists. If $N\ge3$ and $p>N/(N-2)$, then global solutions exist for suitably small data and forcing. The critical exponent $N/(N-2)$ is the Serrin exponent of the elliptic inequality $-\Delta v\ge v^p$, which has no positive solution in $\R^N$ precisely when $p\le N/(N-2)$ \cite{Gid80, SZ02, MP}. Its appearance reflects the fact that, under a positive forcing, global solutions tend to stationary supersolutions, and such supersolutions do not exist below this exponent. The phenomenon has since been extended to time-dependent forcing \cite{JKS, Maj21, AAAAM23, YM}, and to many other operators and geometries \cite{KT26, BBM25, MT26, KM26}. 
It is tempting to replace $|u|^p$ by the gradient nonlinearity $|\nabla u|^q$. The elliptic counterpart then becomes $-\Delta v\ge|\nabla v|^q+\mathbf w$, and an elementary test-function computation (Proposition~\ref{prop:elliptic} below) shows that this inequality has no solution in $\R^N$ when $q\le N/(N-1)$. By analogy with \cite{BLZ}, one might expect that \eqref{eq:main} has no global solution for such $q$. Nonexistence results of this type are indeed known for the combined problem $u_t=\Delta u+|u|^p+|\nabla u|^q+\mathbf w$, where the source term $|u|^p$ is used to control the terms involving $u$ itself \cite{JSS20, Tor23}. For \eqref{eq:main}, however, {the analogy is misleading}. The maximum principle bounds $u$ by $\sup u_0+t\|\mathbf w\|_{L^\infty}$ on every finite time interval, and the same argument applied to difference quotients bounds $\nabla u$. We deduce that \eqref{eq:main} is globally well posed for regular data, whatever $q$ and $\mathbf w$ are (Proposition~\ref{prop:global}). The absence of stationary solutions must therefore show up in a different way, namely through {grow-up}, that is, unboundedness of $u(t)$ as $t\to\infty$. The purpose of this paper is to make this precise.

\begin{defi}\label{def:growup}
A global solution $u$ of \eqref{eq:main} is {bounded} if $\sup_{t\ge0}\|u(t)\|_{L^\infty(\R^N)}<\infty$. It {grows up} if $\sup_{t\ge0}\|u(t)\|_{L^\infty(\R^N)}=\infty$. We say that $u$ {grows linearly on a ball} if there exist $\rho,c>0$ and $C\ge0$ such that
\begin{equation*}
\max_{|x|\le\rho}u(t,x)\ge ct-C\qquad\text{for all }t\ge0.
\end{equation*}
\end{defi}

\subsection{Heuristics and the critical exponent}\label{ss:heuristics}
Two heuristic computations single out the exponent
\begin{equation}\label{eq:qc}
q_c:=\frac{N}{N-1}\quad(N\ge2),\qquad q_c:=+\infty\quad(N=1).
\end{equation}
The first one is scaling. The equation $u_t=\Delta u+|\nabla u|^q$ is invariant under the transformation $u_\lambda(t,x)=\lambda^{\frac{2-q}{q-1}}u(\lambda^2t,\lambda x)$, and the forcing $\mathbf w$ is then transformed into $\lambda^{\frac{2-q}{q-1}+2}\mathbf w(\lambda x)=\lambda^{q'}\mathbf w(\lambda x)$, where $q'=q/(q-1)$. The $L^1$ norm, which carries the relevant quantity $\int\mathbf w$, is multiplied by $\lambda^{q'-N}$. It is therefore invariant exactly when $q'=N$, that is, when $q=q_c$. For $q<q_c$ the forcing becomes stronger at large scales, and for $q>q_c$ it becomes weaker. The same $L^1$-criticality principle produces the Fujita exponent $1+2/N$ and the Bandle--Levine--Zhang exponent $N/(N-2)$. The second heuristic is elliptic. By Proposition~\ref{prop:elliptic}, the stationary problem has no solution when $N\ge2$ and $q\le q_c$. For $q>q_c$ and $N\ge3$, explicit bounded supersolutions exist (Section~\ref{sec:bounded}). Since a bounded global solution should converge to a stationary one, one expects grow-up for $q\le q_c$, and bounded solutions for small data when $q>q_c$. We stress that these computations are only heuristic. The statements below are proved independently of them.

\subsection{Main results}
All results concern the classical solutions provided by Proposition~\ref{prop:global}, with $u_0\in W^{1,\infty}(\R^N)$ and $\mathbf w$ as in \eqref{eq:mass}. Unless stated otherwise, $\mathtt a\equiv1$.

\begin{enumerate}
\item[(A)] \emph{Subcritical range} (Theorem~\ref{thm:sub}). If $1<q<q_c$, every solution grows linearly on a ball. Since $u(t)\le\sup u_0+t\|\mathbf w_+\|_{L^\infty}$, linear growth is the exact order.
\item[(B)] \emph{Critical exponent} (Theorem~\ref{thm:crit}). If $N\ge2$ and $q=q_c$, there exist $\rho,c>0$ and $C\ge0$ such that
\begin{equation*}
\max_{|x|\le\rho}u(t,x)\ge ct-C,\qquad t\ge0.
\end{equation*}
In particular, every solution grows linearly on a ball, and $\liminf_{t\to\infty}t^{-1}\max_{|x|\le\rho}u(t,x)>0$.
\item[(C)] \emph{The plane} (Theorems~\ref{thm:plane} and \ref{thm:plane-local}). If $N=2$, then for {every} $q>1$ every solution grows up. More precisely, $u(t,x)\to+\infty$ locally uniformly, and
\begin{equation*}
\liminf_{t\to\infty}\frac{1}{\log t}\inf_{|x|\le R}u(t,x)\ge\frac{M_\infty}{4\pi}\qquad\text{for every }R>0.
\end{equation*}
\item[(D)] \emph{The quadratic case} (Theorem~\ref{thm:quadratic}). If $q=2$, then $\max_{|x|\le R}u(t,x)\ge\inf u_0-\lambda_Rt$ for all $R,t>0$, where $\lambda_R$ is the principal Dirichlet eigenvalue of $-\Delta-\mathbf w$ in $B_R$. Consequently, $\liminf_{t\to\infty}t^{-1}\sup u(t)\ge-\lambda_0$, where $\lambda_0=\inf\sigma(-\Delta-\mathbf w)$. Moreover $\lambda_0<0$ whenever $N\le2$, $\int\mathbf w\ge0$ and $\mathbf w\not\equiv0$. Combining this with known results on Schr\"odinger semigroups, we obtain that the growth rate is exactly $-\lambda_0$ (Remark~\ref{rem:exact}).
\item[(E)] \emph{Supercritical range} (Theorem~\ref{thm:bounded}). If $N\ge3$ and $q>q_c$, then for small nonnegative $u_0$ and $\mathbf w$, possibly with $\int\mathbf w>0$, the solution satisfies $0\le u(t,x)\le\delta\jb{x}^{-\gamma}$ for all $t\ge0$. If moreover $\mathtt a\equiv1$ and $u_0\equiv0$, it increases monotonically to a bounded limit. Conversely, for every $q>1$ and every $N$, a scale-invariant largeness condition on $\mathbf w$ forces linear growth (Theorem~\ref{thm:large}).
\item[(F)] \emph{Time-dependent forcing} (Theorem~\ref{thm:profiles}). For every $q>1$ and every $N$, if $A(t)=\int_0^t\mathtt a$ satisfies $A(t)/t\to\infty$, then $\max_{|x|\le\rho}u(t,x)\ge cA(t)$ for large $t$. If $1<q\le q_c$ and $\liminf_{t\to\infty}A(t)/t>0$, solutions grow linearly on a ball. For $\mathtt a(t)\sim t^\sigma$ with $-1<\sigma<0$, we prove grow-up with a power rate when $q<N/(N-1-\sigma)$ (Corollary~\ref{cor:decay}). This is a partial result, since scaling suggests the threshold $(N-2\sigma)/(N-1-2\sigma)$ (Conjecture~\ref{conj:decay}).
\end{enumerate}

Table~\ref{tab:summary} in Section~\ref{sec:open} summarises the situation and indicates what remains open. In the supercritical range, the outcome depends on the size of the data, just as in the Fujita problem. For $q=2$, the growth rate is governed by $\lambda_0$, and Hardy's inequality gives a sufficient condition for $\lambda_0\ge0$ in terms of the size of $\mathbf w$ relative to $(N-2)^2/(4|x|^2)$ (Corollary~\ref{cor:quadratic}).

\subsection{Relation to the literature}
The Cauchy problem for $u_t=\Delta u+|\nabla u|^q$ and its variants has been studied in detail. The local theory, including rough data, is developed in \cite{AB98, BASW02, BL99, GGK03}, and the large-time behaviour of solutions is described in \cite{Gil05}. Gradient blow-up, which requires boundary conditions, is surveyed in \cite{Sou01}. The growth of the mass $\int u$ was studied in \cite{LS03, LS05}, and the equation with gradient absorption in \cite{BAK99}. Souplet and Zhang \cite{SZ06} considered the inhomogeneous equation $u_t-\Delta u=|\nabla u|^p+h(x)$, and related the existence of global, bounded or locally bounded solutions to the existence of stationary solutions. Competition between a reaction term and a gradient term, as in the Chipot--Weissler equation $u_t=\Delta u+|u|^p-\mu|\nabla u|^q$, was studied in \cite{CW89, KP89, Fil91, STW96, SW97, SW99, STW99}. The influence of a gradient source on the Fujita exponent is analysed in \cite{JSS20, Tor23}. Such equations fit into the theory of viscosity solutions \cite{CIL92}.

Viscous Hamilton--Jacobi equations are usually written in the form $U_t-\Delta U+|\nabla U|^q=f(x)$. In the notation of \eqref{eq:main}, we have $U=-u$ and $f=-\mathbf w$. The large-time behaviour $U(t)-\lambda t\to\phi$, governed by an ergodic problem, has been established in periodic settings and in bounded domains \cite{BS01, BPT10, Tch10}, and for coercive potentials $f(x)\to+\infty$ in $\R^N$ \cite{KS06, Ich12, Ich15, BM16, BQR21, QR24}. Uniqueness questions for the ergodic problem are treated in \cite{ABC19}. The present problem corresponds to a potential $f=-\mathbf w$ that {decays} at infinity, which is the non-coercive regime. In this regime, recurrence and transience of the underlying controlled diffusion become decisive, and so does the dimension. We refer to the work of Ichihara on criticality \cite{Ich13}, and to \cite{Gri99} for the role of recurrence and capacity. For $q=2$, this is the classical dichotomy between critical and subcritical Schr\"odinger operators \cite{Mur86, Pin07}. It is also the weak coupling phenomenon: in dimensions one and two, every potential with positive mass produces a bound state \cite{Sim76}.

To the best of our knowledge, the Fujita-type dichotomy for \eqref{eq:main}, with the critical exponent $q_c=N/(N-1)$ separating grow-up from boundedness, has not been formulated before. The exponent $N/(N-1)$ appears in \cite{JSS20, Tor23} as a blow-up threshold for the combined problem with the source term $|u|^p$, and the ergodic problem with potentials that vanish at infinity is studied in \cite{Ich13, CI16}. Here we show that, without the source $|u|^p$, finite-time blow-up is replaced by grow-up.

Our contribution is to organise these phenomena around the Fujita-type critical exponent $q_c=N/(N-1)$, and to prove the growth and boundedness statements (A)--(F) by elementary means. The only tools are test functions, the comparison principle and principal eigenvalues on balls.

\subsection{Organisation}
Section~\ref{sec:prelim} collects preliminaries, namely cut-off functions, the comparison principle, global well-posedness, the stationary problem and the Cole--Hopf transformation. Sections~\ref{sec:sub}, \ref{sec:crit}, \ref{sec:plane}, \ref{sec:quadratic} and \ref{sec:bounded} contain the proofs of (A)--(F). Section~\ref{sec:open} contains a summary, some comments and open problems.

\subsection{Notation}
$B_R$ is the open ball of radius $R$ centred at $0$, $A_R:=\{R\le|x|\le2R\}$, $\jb{x}:=(1+|x|^2)^{1/2}$, $q':=q/(q-1)$, and $\mathbf w_\pm:=\max\{\pm\mathbf w,0\}$. We write $\|\cdot\|_\infty$ for the norm of $L^\infty(\R^N)$, and $|\mS^{N-1}|$ for the surface measure of the unit sphere. For $R>0$, we set $M(R):=\int\mathbf w\,\varphi_R\,\dd x$, where $\varphi_R$ is the cut-off function of Lemma~\ref{lem:cutoff}. By dominated convergence, $M(R)\to M_\infty$ as $R\to\infty$.

\section{Preliminaries}\label{sec:prelim}

\subsection{Cut-off functions}
We fix a radial function $\varphi_0\in C^\infty(\R^N)$, nonincreasing in $|x|$, with $0\le\varphi_0\le1$, $\varphi_0=1$ on $B_1$ and $\varphi_0=0$ outside $B_2$. Given $q>1$, we fix an integer $k\ge\max\{2,q'\}$, and we set $\varphi:=\varphi_0^k$ and $\varphi_R(x):=\varphi(x/R)$. Since $\varphi_0$ is radially nonincreasing, the map $R\mapsto\varphi_R(x)$ is nondecreasing, and $\varphi_R\to1$ pointwise as $R\to\infty$.

\begin{lem}\label{lem:cutoff}
There are constants $\Theta,C_*>0$, depending only on $N$, $q$ and $\varphi_0$, such that, for every $R>0$, the following hold.
\begin{enumerate}
\item[(a)] $\int_{\R^N}|\nabla\varphi_R|^{q'}\varphi_R^{-q'/q}\,\dd x=\Theta\,R^{N-q'}$, where the integrand is set equal to $0$ where $\varphi_R=0$.
\item[(b)] $|\nabla\varphi_R|\le C_*R^{-1}\varphi_R^{1/q}\,\mathbf 1_{A_R}$.
\item[(c)] $|\Delta\varphi_R|\le C_*R^{-2}\,\mathbf 1_{A_R}$.
\end{enumerate}
\end{lem}

\begin{proof}
We have $\nabla\varphi=k\varphi_0^{k-1}\nabla\varphi_0$, which is supported in $A_1$. Since $k\ge q'$, we have $(k-1)/k\ge1/q$, and since $0\le\varphi\le1$, we get $\varphi_0^{k-1}=\varphi^{(k-1)/k}\le\varphi^{1/q}$. Hence $|\nabla\varphi|\le C\varphi^{1/q}\mathbf 1_{A_1}$, and it follows that $|\nabla\varphi|^{q'}\varphi^{-q'/q}\le C^{q'}\mathbf 1_{A_1}$ is bounded and compactly supported. In the same way, $\Delta\varphi=k\varphi_0^{k-1}\Delta\varphi_0+k(k-1)\varphi_0^{k-2}|\nabla\varphi_0|^2$, so that $|\Delta\varphi|\le C\mathbf 1_{A_1}$ because $k\ge2$. The constant $\Theta=\int|\nabla\varphi|^{q'}\varphi^{-q'/q}$ is positive since $\varphi$ is not constant. The three statements follow by scaling, since $\nabla\varphi_R=R^{-1}(\nabla\varphi)(\cdot/R)$ and $\Delta\varphi_R=R^{-2}(\Delta\varphi)(\cdot/R)$.
\end{proof}

We note for later use that a nonnegative $C^1$ function has a vanishing gradient at each of its zeros. In particular, $\nabla\varphi_R=0$ wherever $\varphi_R=0$.

In the plane, we also use logarithmic cut-off functions. They exploit the parabolicity of $\R^2$, namely the fact that the $2$-capacity of a ball relative to infinity vanishes \cite{Gri99}.

\begin{lem}\label{lem:log}
Let $N=2$, and let $\eta\in C^\infty(\R)$ be nonincreasing with $\eta=1$ on $(-\infty,1]$ and $\eta=0$ on $[2,\infty)$. For $R\ge2$, define $\ell_R(x):=\eta(\log|x|/\log R)$ for $|x|\ge1$ and $\ell_R(x):=1$ for $|x|<1$. Then $\ell_R\in C_c^\infty(\R^2)$, $0\le\ell_R\le1$, $\ell_R=1$ on $B_R$, $\operatorname{supp}\ell_R\subset\overline{B_{R^2}}$, $\ell_R\to1$ pointwise as $R\to\infty$, and
\begin{equation*}
\int_{\R^2}|\Delta\ell_R|\,\dd x=\frac{2\pi\|\eta''\|_{L^1}}{\log R},\qquad \int_{\R^2}|\nabla\ell_R|^2\,\dd x=\frac{2\pi\|\eta'\|_{L^2}^2}{\log R}.
\end{equation*}
\end{lem}

\begin{proof}
Write $L=\log R$ and $r=|x|$. For $r\ge1$, we have $\ell_R=f(r)$ with $f(r)=\eta(\log r/L)$, and $f(r)=1$ for $1\le r\le R$. This gives smoothness, the support property and the pointwise convergence. Setting $s=\log r/L$, we have $f'(r)=\eta'(s)/(rL)$ and $f''(r)=\eta''(s)/(r^2L^2)-\eta'(s)/(r^2L)$. In two dimensions, we obtain
\begin{equation*}
\Delta\ell_R=f''+\frac{f'}{r}=\frac{\eta''(s)}{r^2L^2},
\end{equation*}
the first-order terms cancelling because $\log|x|$ is harmonic in the plane. Using polar coordinates and $\dd r/r=L\,\dd s$, we get $\int|\Delta\ell_R|=2\pi\int_R^{R^2}|\eta''(s)|\,\dd r/(rL^2)=2\pi\|\eta''\|_{L^1}/L$. The gradient integral is computed in the same way from $|\nabla\ell_R|=|\eta'(s)|/(rL)$.
\end{proof}

At the critical exponent, we need a logarithmic cut-off function adapted to the weight appearing in Young's inequality.

\begin{lem}\label{logcutoff}
Let $N\ge2$ and $q=q_c=N/(N-1)$, so that $q'=N$. Let $\eta$ be as in Lemma~\ref{lem:log}. For every $R>1$, define
\begin{equation*}
\zeta_R(x)=
\begin{cases}
1, & |x|\le 1,\\[1mm]
\Big(\eta\Big(\dfrac{\log|x|}{\log R}\Big)\Big)^N, & |x|>1.
\end{cases}
\end{equation*}
Then the following properties hold.
\begin{enumerate}
\item[(a)] $\zeta_R\in C_c^\infty(\R^N)$, $0\le\zeta_R\le1$, $\zeta_R=1$ in $B_R$, and $\operatorname{supp}\zeta_R\subset\overline{B_{R^2}}$. Moreover, $\zeta_R(x)\to1$ as $R\to\infty$ for every $x\in\R^N$.
\item[(b)] For every $x\in\R^N$,
\begin{equation*}
|\nabla\zeta_R(x)|\le\frac{N\|\eta'\|_\infty}{|x|\log R}\,\zeta_R(x)^{1/q}\,\mathbf 1_{\{R\le|x|\le R^2\}}(x).
\end{equation*}
\item[(c)] With the integrand set equal to $0$ where $\zeta_R=0$, we have
\begin{equation*}
\int_{\R^N}|\nabla\zeta_R|^{q'}\zeta_R^{-q'/q}\,\dd x=\frac{\Theta_N}{(\log R)^{N-1}},\qquad \Theta_N:=N^N|\mS^{N-1}|\int_1^2|\eta'(s)|^N\,\dd s>0.
\end{equation*}
In particular, this integral tends to $0$ as $R\to\infty$.
\end{enumerate}
\end{lem}

\begin{proof}
Since $\eta=1$ on $(-\infty,1]$ and $\eta=0$ on $[2,\infty)$, we have $\zeta_R(x)=1$ for $|x|\le R$ and $\zeta_R(x)=0$ for $|x|\ge R^2$. In particular, the two expressions in the definition agree near $|x|=1$, and we infer that $\zeta_R\in C_c^\infty(\R^N)$. For fixed $x$, we have $\zeta_R(x)=1$ as soon as $R>|x|$, which proves (a).

The gradient of $\zeta_R$ vanishes outside the annulus $\{R<|x|<R^2\}$. Inside it, writing $s=\log|x|/\log R$, the chain rule gives
\begin{equation*}
|\nabla\zeta_R(x)|=\frac{N}{|x|\log R}\,\eta(s)^{N-1}|\eta'(s)|.
\end{equation*}
Since $q=N/(N-1)$, we have $\zeta_R^{1/q}=\eta(s)^{N-1}$, and (b) follows from $|\eta'|\le\|\eta'\|_\infty$.

For (c), we use $q'=N$ and $q'/q=N-1$. Wherever $\zeta_R>0$, the previous formula yields
\begin{equation*}
|\nabla\zeta_R(x)|^{q'}\zeta_R(x)^{-q'/q}=\frac{N^N}{|x|^N(\log R)^N}\,|\eta'(s)|^N .
\end{equation*}
This identity also holds, with the zero convention, where $\zeta_R=0$, because $\eta'(s)=0$ there. Using polar coordinates and the change of variables $s=\log r/\log R$, for which $\dd r/r=(\log R)\,\dd s$, we obtain
\begin{equation*}
\int_{\R^N}|\nabla\zeta_R|^{q'}\zeta_R^{-q'/q}\,\dd x=\frac{N^N|\mS^{N-1}|}{(\log R)^N}\int_R^{R^2}\Big|\eta'\Big(\frac{\log r}{\log R}\Big)\Big|^N\frac{\dd r}{r}=\frac{N^N|\mS^{N-1}|}{(\log R)^{N-1}}\int_1^2|\eta'(s)|^N\,\dd s .
\end{equation*}
The constant $\Theta_N$ is finite, and it is positive because $\eta$ decreases from $1$ to $0$ on $[1,2]$. Since $N\ge2$, the right-hand side tends to $0$ as $R\to\infty$. This completes the proof.
\end{proof}

\subsection{Maximum and comparison principles}
We use the following maximum principle for functions that are bounded from above on the whole space. The proof is standard, see \cite{Fri64, Lie96, QS}. We include it because the drift coefficients used below are only bounded and measurable.

\begin{lem}\label{lem:maxprinciple}
Let $T>0$ and let $b:(0,T]\times\R^N\to\R^N$ be measurable with $|b|\le B$. Let $z\in C([0,T]\times\R^N)\cap C^{1,2}((0,T]\times\R^N)$ be bounded from above, and suppose that $z_t-\Delta z-b\cdot\nabla z\le0$ pointwise in $(0,T]\times\R^N$ and $z(0,\cdot)\le0$. Then $z\le0$ in $[0,T]\times\R^N$.
\end{lem}

\begin{proof}
Let $\kappa>2N+B$ and $P(t,x)=(1+|x|^2)e^{\kappa t}$. Since $2|b\cdot x|\le B(1+|x|^2)$, we obtain
\begin{equation*}
P_t-\Delta P-b\cdot\nabla P\ge e^{\kappa t}\big((\kappa-B)(1+|x|^2)-2N\big)>0 .
\end{equation*}
Fix $\varepsilon>0$ and set $Z=z-\varepsilon P$. Since $z$ is bounded from above, there is $\rho>0$ such that $Z<0$ on $[0,T]\times\{|x|\ge\rho\}$. Moreover $Z(0,\cdot)\le-\varepsilon$. Suppose that $Z$ takes a nonnegative value. The function $m(t)=\max_{|x|\le\rho}Z(t,x)$ is continuous with $m(0)<0$. Hence there is a first time $t_0\in(0,T]$ with $m(t_0)=0$, attained at some $x_0$ with $|x_0|<\rho$. Since $Z(t_0,\cdot)\le0$ on $\R^N$, the point $x_0$ is a global spatial maximum, and we get $\nabla Z(t_0,x_0)=0$ and $\Delta Z(t_0,x_0)\le0$. Since $Z(t,x_0)<0$ for $t<t_0$, we also have $Z_t(t_0,x_0)\ge0$. It follows that $Z_t-\Delta Z-b\cdot\nabla Z\ge0$ at $(t_0,x_0)$. On the other hand, by construction this quantity is at most $-\varepsilon(P_t-\Delta P-b\cdot\nabla P)<0$. This contradiction shows that $z\le\varepsilon P$ for every $\varepsilon>0$, and we conclude that $z\le0$.
\end{proof}

\begin{lem}\label{lem:comparison}
Let $T>0$ and $g\in C([0,T]\times\R^N)$. Let $u_1,u_2\in C([0,T]\times\R^N)\cap C^{1,2}((0,T]\times\R^N)$ be bounded, with bounded gradients on $(0,T]\times\R^N$. Suppose that
\begin{equation*}
\partial_tu_1-\Delta u_1-|\nabla u_1|^q-g\le0\le\partial_tu_2-\Delta u_2-|\nabla u_2|^q-g\quad\text{in }(0,T]\times\R^N,
\end{equation*}
and $u_1(0)\le u_2(0)$. Then $u_1\le u_2$ in $[0,T]\times\R^N$.
\end{lem}

\begin{proof}
Let $z=u_1-u_2$. Since $F(\xi)=|\xi|^q$ is of class $C^1$ with $\nabla F(\xi)=q|\xi|^{q-2}\xi$, we have $|\nabla u_1|^q-|\nabla u_2|^q=b\cdot\nabla z$, where $b=\int_0^1\nabla F(\theta\nabla u_1+(1-\theta)\nabla u_2)\,\dd\theta$. If $L$ bounds both gradients, then $|b|\le qL^{q-1}$. Thus $z_t-\Delta z-b\cdot\nabla z\le0$ and $z(0)\le0$, and Lemma~\ref{lem:maxprinciple} gives $z\le0$.
\end{proof}

\subsection{Global well-posedness}
The following proposition shows that finite-time blow-up is excluded for \eqref{eq:main}. It is the starting point of the whole analysis.

\begin{prop}\label{prop:global}
Let $q>1$, $u_0\in W^{1,\infty}(\R^N)$, $\mathbf w\in W^{1,\infty}(\R^N)$ and $\mathtt a\in C([0,\infty))$. Then \eqref{eq:main} has a unique global solution $u\in C([0,\infty)\times\R^N)\cap C^{1,2}((0,\infty)\times\R^N)$ such that $u$ and $\nabla u$ are bounded on $(0,T]\times\R^N$ for every $T>0$. Moreover, for all $t\ge0$ and $x\in\R^N$,
\begin{equation}\label{eq:apriori}
\begin{gathered}
\inf u_0-\|\mathbf w\|_\infty\int_0^t|\mathtt a|\le u(t,x)\le\sup u_0+\|\mathbf w\|_\infty\int_0^t|\mathtt a|,\\
\|\nabla u(t)\|_\infty\le\|\nabla u_0\|_\infty+\|\nabla\mathbf w\|_\infty\int_0^t|\mathtt a| .
\end{gathered}
\end{equation}
When $\mathtt a\equiv1$, the first estimate can be sharpened to $\inf u_0-t\|\mathbf w_-\|_\infty\le u\le\sup u_0+t\|\mathbf w_+\|_\infty$. In particular, $u\ge\inf u_0$ if moreover $\mathbf w\ge0$.
\end{prop}

\begin{proof}
Local existence and uniqueness are standard. One writes \eqref{eq:main} in the Duhamel form
\begin{equation*}
u(t)=e^{t\Delta}u_0+\int_0^te^{(t-s)\Delta}\big(|\nabla u(s)|^q+\mathtt a(s)\mathbf w\big)\,\dd s .
\end{equation*}
Since $\xi\mapsto|\xi|^q$ is locally Lipschitz and $\|\nabla e^{t\Delta}f\|_\infty\le Ct^{-1/2}\|f\|_\infty$, a contraction argument in $L^\infty((0,T);W^{1,\infty}(\R^N))$ gives a solution on an interval whose length depends only on $\|u_0\|_{W^{1,\infty}}$ and on bounds for $\mathtt a$ and $\mathbf w$. The forcing is continuous in $t$ and Lipschitz in $x$, so that parabolic regularity shows that this solution is classical. We refer to \cite{AB98, BASW02, GGK03} for details in the unforced case. The bounded Lipschitz forcing term does not affect the arguments. In particular, the maximal solution is global as soon as $\|u(t)\|_{W^{1,\infty}}$ stays bounded on bounded time intervals. Uniqueness in the stated class follows from Lemma~\ref{lem:comparison}.

It remains to prove \eqref{eq:apriori} on $[0,T]$ for each $T$ smaller than the maximal existence time. The functions $\overline u(t)=\sup u_0+\|\mathbf w\|_\infty\int_0^t|\mathtt a|$ and $\underline u(t)=\inf u_0-\|\mathbf w\|_\infty\int_0^t|\mathtt a|$ have zero gradient, and they are, respectively, a supersolution and a subsolution of \eqref{eq:main}. Lemma~\ref{lem:comparison} with $g=\mathtt a\mathbf w$ gives the bounds on $u$. When $\mathtt a\equiv1$, the sharper bounds follow in the same way from the supersolution $\sup u_0+t\|\mathbf w_+\|_\infty$ and the subsolution $\inf u_0-t\|\mathbf w_-\|_\infty$, since $|\nabla u|^q\ge0$.

For the gradient bound, fix a unit vector $e$ and $h\ne0$, and let $v_h(t,x)=(u(t,x+he)-u(t,x))/h$. Arguing as in the proof of Lemma~\ref{lem:comparison}, we see that $v_h$ satisfies
\begin{equation*}
\partial_tv_h-\Delta v_h-b_h\cdot\nabla v_h=\mathtt a(t)\,\frac{\mathbf w(x+he)-\mathbf w(x)}{h},
\end{equation*}
with a bounded drift $b_h$. The right-hand side is bounded in absolute value by $|\mathtt a(t)|\,\|\nabla\mathbf w\|_\infty$. Hence the functions $z_\pm=\pm v_h-\|\nabla u_0\|_\infty-\|\nabla\mathbf w\|_\infty\int_0^t|\mathtt a|$ are bounded from above and satisfy $\partial_tz_\pm-\Delta z_\pm-b_h\cdot\nabla z_\pm\le0$ and $z_\pm(0)\le0$. Lemma~\ref{lem:maxprinciple} gives $|v_h|\le\|\nabla u_0\|_\infty+\|\nabla\mathbf w\|_\infty\int_0^t|\mathtt a|$, and letting $h\to0$ we get the gradient bound. Since the right-hand sides of \eqref{eq:apriori} are finite for every $t$, we conclude that the solution is global.
\end{proof}

\begin{rem}
Proposition~\ref{prop:global} holds for {every} $q>1$ and {every} bounded Lipschitz forcing, whatever its sign or mass. This is the fundamental difference with the semilinear problem $u_t=\Delta u+|u|^p+\mathbf w$ of \cite{BLZ}. The nonlinearity $|\nabla u|^q$ vanishes on constants, and it cannot by itself increase the maximum of $u$. Only the forcing can do so. Consequently, the only way for a positive forcing to destabilise the problem is to make $u$ unbounded as $t\to\infty$.
\end{rem}

\subsection{The stationary problem}
The following Liouville-type statement is the elliptic counterpart of our results. Its proof contains, in the simplest setting, the test-function mechanism used later.

\begin{prop}\label{prop:elliptic}
Let $q>1$ and $\mathbf w\in L^1(\R^N)\cap C(\R^N)$ with $\int\mathbf w>0$.
\begin{enumerate}
\item[(i)] If $N=1$, or if $N\ge2$ and $q\le q_c=N/(N-1)$, there is no $v\in C^2(\R^N)$ such that $-\Delta v\ge|\nabla v|^q+\mathbf w$ in $\R^N$. No sign or growth condition on $v$ is assumed.
\item[(ii)] If $N=2$ and $\mathbf w\ge0$, then for every $q>1$ there is no $v\in C^2(\R^2)$, bounded from below, such that $-\Delta v\ge|\nabla v|^q+\mathbf w$ in $\R^2$.
\end{enumerate}
\end{prop}

\begin{proof}
(i) We multiply the inequality by $\varphi_R$ and integrate by parts once. On the set where $\varphi_R>0$, Young's inequality gives
\begin{equation*}
|\nabla v||\nabla\varphi_R|=\big(|\nabla v|\varphi_R^{1/q}\big)\big(|\nabla\varphi_R|\varphi_R^{-1/q}\big)\le\frac12|\nabla v|^q\varphi_R+C_q|\nabla\varphi_R|^{q'}\varphi_R^{-q'/q},
\end{equation*}
and both sides vanish where $\varphi_R=0$. Using Lemma~\ref{lem:cutoff}(a), we obtain
\begin{equation}\label{eq:elliptic1}
\frac12\int|\nabla v|^q\varphi_R+M(R)\le C_q\Theta R^{N-q'} .
\end{equation}
If $q'>N$, which always holds when $N=1$ and is equivalent to $q<q_c$ when $N\ge2$, the right-hand side tends to zero as $R\to\infty$, while $M(R)\to M_\infty>0$. This is a contradiction. If $q'=N$, then \eqref{eq:elliptic1} gives $\int|\nabla v|^q\varphi_R\le2C_q\Theta$ for $R$ large, since $M(R)>0$ for such $R$. By monotone convergence, we get $\int_{\R^N}|\nabla v|^q<\infty$. Instead of Young's inequality, we now use Lemma~\ref{lem:cutoff}(b) and H\"older's inequality to get
\begin{equation*}
\int|\nabla v||\nabla\varphi_R|\le C_*R^{-1}|A_R|^{1/q'}\Big(\int_{A_R}|\nabla v|^q\Big)^{1/q}=C\Big(\int_{A_R}|\nabla v|^q\Big)^{1/q},
\end{equation*}
because $R^{-1}|A_R|^{1/q'}$ does not depend on $R$ when $q'=N$. The right-hand side tends to $0$ as $R\to\infty$, being bounded by a power of the tail of a convergent integral. Since $M(R)\le\int\nabla v\cdot\nabla\varphi_R-\int|\nabla v|^q\varphi_R\le\int|\nabla v||\nabla\varphi_R|$, we deduce that $M_\infty\le0$, again a contradiction.

(ii) If $\mathbf w\ge0$, then $v$ is superharmonic in $\R^2$ and bounded from below. Hence $v$ is constant, see for example \cite{AG01}. The inequality then gives $0\ge\mathbf w$, which contradicts $\int\mathbf w>0$.
\end{proof}

\begin{rem}
Proposition~\ref{prop:elliptic}(i) is universal, in the sense that it requires no information on $v$. In the parabolic problem, the corresponding computation produces the extra term $\frac{\dd}{\dd t}\int u\varphi_R$, which cannot be absorbed by the gradient term. The main observation of this paper is that this term need not be absorbed at all. It measures the growth of $u$, and keeping it yields the grow-up results of Sections~\ref{sec:sub} and \ref{sec:crit}. For $q>q_c$ and $N\ge3$, bounded positive supersolutions do exist, see Lemma~\ref{lem:supersolution}. With $\mathbf w\equiv0$, the argument above shows that every $C^2$ solution of $-\Delta v\ge|\nabla v|^q$ in $\R^N$ is constant when $q\le q_c$, see also \cite{MP}. For the equation $-\Delta v=|\nabla v|^q$, Lions proved that every solution in $\R^N$ is constant for every $q>1$ \cite{Lio85}, see also \cite{BGV14}.
\end{rem}

\subsection{The Cole--Hopf transformation and principal eigenvalues}\label{ss:colehopf}
When $q=2$, the transformation $v=e^u$ \cite{Hop50, Col51} turns \eqref{eq:main} into a linear equation. Indeed, if $u$ is the solution given by Proposition~\ref{prop:global}, then $v_t=e^uu_t$ and $\Delta v=e^u(\Delta u+|\nabla u|^2)$, so that
\begin{equation}\label{eq:colehopf}
v_t=\Delta v+\mathtt a(t)\mathbf w(x)\,v,\qquad v>0 .
\end{equation}
The growth of $u$ is thus governed by the spectral properties of the Schr\"odinger operator $-\Delta-\mathbf w$. For $R>0$, let
\begin{equation*}
\lambda_R:=\min\Big\{\int_{B_R}\big(|\nabla f|^2-\mathbf w f^2\big)\,\dd x:\ f\in H^1_0(B_R),\ \|f\|_{L^2(B_R)}=1\Big\}
\end{equation*}
be the principal Dirichlet eigenvalue of $-\Delta-\mathbf w$ in $B_R$. Since $\mathbf w\in W^{1,\infty}$, the minimum is attained by an eigenfunction $\psi_R\in C^{2}(\overline{B_R})$ with $\psi_R>0$ in $B_R$ and $\psi_R=0$ on $\partial B_R$, see for instance \cite[Section 6.5]{Eva10} and \cite{GT01}. In particular, $\partial_\nu\psi_R\le0$ on $\partial B_R$. We normalise $\psi_R$ by $\int_{B_R}\psi_R=1$. Since $H^1_0(B_R)\subset H^1_0(B_{R'})$ for $R<R'$, the map $R\mapsto\lambda_R$ is nonincreasing, and $\lambda_R\ge-\|\mathbf w\|_\infty$. By density of $C_c^\infty(\R^N)$ in $H^1(\R^N)$, we have
\begin{equation}\label{eq:lambda0}
\lambda_0:=\lim_{R\to\infty}\lambda_R=\inf\Big\{\int\big(|\nabla f|^2-\mathbf wf^2\big):\ f\in C_c^\infty(\R^N),\ \|f\|_{L^2}=1\Big\},
\end{equation}
which is the bottom of the spectrum of the self-adjoint operator $-\Delta-\mathbf w$ on $L^2(\R^N)$ \cite{RS78}.

\section{Linear growth in the subcritical range}\label{sec:sub}

Throughout Sections~\ref{sec:sub}--\ref{sec:quadratic}, $u$ is the global solution of \eqref{eq:main} given by Proposition~\ref{prop:global}, with $u_0\in W^{1,\infty}(\R^N)$, $\mathbf w$ satisfying \eqref{eq:mass}, and $\mathtt a\in C([0,\infty))$. We set $A(t):=\int_0^t\mathtt a(s)\,\dd s$,
\begin{equation*}
F_R(t):=\int_{\R^N}u(t,x)\,\varphi_R(x)\,\dd x,\qquad C_q:=\frac{q-1}{q}\Big(\frac2q\Big)^{\frac1{q-1}} .
\end{equation*}
The constant $C_q$ is chosen so that Young's inequality takes the form $xy\le\frac12x^q+C_qy^{q'}$ for $x,y\ge0$. We shall often use the following elementary observation. If $\psi\ge0$ is continuous with support in $\overline{B_\rho}$, then
\begin{equation}\label{eq:avg}
\int u(t)\psi\le\Big(\int\psi\Big)\max_{|x|\le\rho}u(t,x),\qquad \int u_0\psi\ge\Big(\int\psi\Big)\inf u_0 ,
\end{equation}
whatever the sign of $u$. The whole analysis rests on the following identities.

\begin{lem}\label{lem:identity}
\begin{enumerate}
\item[(a)] For every $R>0$, we have $F_R\in C([0,\infty))\cap C^1((0,\infty))$, and for $t>0$,
\begin{align}
F_R'(t)&=\int u\,\Delta\varphi_R+\int|\nabla u|^q\varphi_R+\mathtt a(t)M(R)\label{eq:id1}\\
&\ge\frac12\int|\nabla u|^q\varphi_R+\mathtt a(t)M(R)-C_q\Theta R^{N-q'} .\label{eq:id2}
\end{align}
\item[(b)] Let $N\ge2$ and $q=q_c$. For $R>1$, let $\zeta_R$ be the function of Lemma~\ref{logcutoff}, and set
\begin{equation*}
G_R(t):=\int_{\R^N}u(t,x)\,\zeta_R(x)\,\dd x,\qquad \widetilde M(R):=\int_{\R^N}\mathbf w(x)\,\zeta_R(x)\,\dd x .
\end{equation*}
Then $\widetilde M(R)\to M_\infty$ as $R\to\infty$, $G_R\in C([0,\infty))\cap C^1((0,\infty))$, and for $t>0$,
\begin{equation}\label{eq:id3}
G_R'(t)\ge\frac12\int|\nabla u|^q\zeta_R+\mathtt a(t)\widetilde M(R)-\frac{C_q\Theta_N}{(\log R)^{N-1}} .
\end{equation}
\end{enumerate}
\end{lem}

\begin{proof}
(a) Since $u_t$ is continuous in $(0,\infty)\times\R^N$ and $\varphi_R$ has compact support, $F_R$ is of class $C^1$ on $(0,\infty)$ with $F_R'=\int u_t\varphi_R$. Continuity at $t=0$ follows from the continuity of $u$ on $[0,\infty)\times\R^N$. Inserting the equation and integrating by parts twice, we obtain \eqref{eq:id1}. Integrating by parts only once, we get instead
\begin{equation*}
F_R'=-\int\nabla u\cdot\nabla\varphi_R+\int|\nabla u|^q\varphi_R+\mathtt a(t)M(R).
\end{equation*}
As in the proof of Proposition~\ref{prop:elliptic}, Young's inequality gives
\begin{equation*}
|\nabla u||\nabla\varphi_R|\le\frac12|\nabla u|^q\varphi_R+C_q|\nabla\varphi_R|^{q'}\varphi_R^{-q'/q},
\end{equation*}
and \eqref{eq:id2} follows from Lemma~\ref{lem:cutoff}(a).

(b) The convergence $\widetilde M(R)\to M_\infty$ follows from dominated convergence, since $0\le\zeta_R\le1$ and $\zeta_R\to1$ pointwise. The regularity of $G_R$ is proved as in (a). Integrating by parts once and using Young's inequality in the same way, we get \eqref{eq:id3} from Lemma~\ref{logcutoff}(c).
\end{proof}

Note that \eqref{eq:id2} is the parabolic analogue of \eqref{eq:elliptic1}, with the time derivative $F_R'$ in place of $0$. We use it for all $q$, and \eqref{eq:id1} in Section~\ref{sec:plane}.

\begin{thm}\label{thm:sub}
Let $\mathtt a\equiv1$ and $1<q<q_c$, where $q_c$ is given by \eqref{eq:qc}. Then there is $R_0>0$ such that, for every $R\ge R_0$ and $t\ge0$,
\begin{equation}\label{eq:linear}
F_R(t)\ge F_R(0)+\frac{M_\infty}{2}\,t .
\end{equation}
Consequently, for $R\ge R_0$ and $t\ge0$,
\begin{equation}\label{eq:linball}
\max_{|x|\le2R}u(t,x)\ge\inf u_0+\frac{M_\infty}{2\|\varphi_R\|_{L^1}}\,t ,
\end{equation}
and in particular $\max_{|x|\le2R}u(t,x)\ge\frac{M_\infty}{4|B_{2R}|}\,t$ for $t\ge t_R:=4\|u_0\|_\infty|B_{2R}|/M_\infty$. Together with $u(t,x)\le\sup u_0+t\|\mathbf w_+\|_\infty$, this shows that $u$ grows exactly linearly on the ball $B_{2R}$.
\end{thm}

\begin{proof}
Since $q<q_c$ means $q'>N$, the term $C_q\Theta R^{N-q'}$ tends to zero as $R\to\infty$, while $M(R)\to M_\infty>0$. We choose $R_0$ such that $M(R)-C_q\Theta R^{N-q'}\ge M_\infty/2$ for $R\ge R_0$. Dropping the nonnegative gradient term in \eqref{eq:id2}, we obtain $F_R'\ge M_\infty/2$ on $(0,\infty)$. Integrating on $(\varepsilon,t)$ and letting $\varepsilon\to0$, we get \eqref{eq:linear}. Applying \eqref{eq:avg} with $\psi=\varphi_R$ and $\rho=2R$, we deduce \eqref{eq:linball}. Finally, $\|\varphi_R\|_{L^1}\le|B_{2R}|$ and $\inf u_0\ge-\|u_0\|_\infty$, so that \eqref{eq:linball} gives
\begin{equation*}
\max_{|x|\le2R}u(t,x)\ge-\|u_0\|_\infty+\frac{M_\infty}{2|B_{2R}|}\,t\ge\frac{M_\infty}{4|B_{2R}|}\,t\qquad\text{for }t\ge t_R .
\end{equation*}
The upper bound is part of Proposition~\ref{prop:global}.
\end{proof}

\begin{rem}
Three features of Theorem~\ref{thm:sub} are worth noting. First, no sign condition is imposed on $\mathbf w$ or on $u_0$. The conclusion holds for sign-changing forcing with positive mass and for arbitrary Lipschitz initial data. Second, the proof uses the gradient nonlinearity only through the inequality $|\nabla u|^q\ge0$ and to absorb the flux term $\int\nabla u\cdot\nabla\varphi_R$. The mass supplied by the forcing is the driving mechanism, and the gradient term cannot evacuate it fast enough when $q'>N$. Third, integrating \eqref{eq:id2} with $\mathtt a\equiv1$ gives, for every $q>1$ and $R>0$,
\begin{equation*}
\frac1t\int_0^t\!\!\int|\nabla u|^q\varphi_R\le\frac2t\big(F_R(t)-F_R(0)\big)+2C_q\Theta R^{N-q'}-2M(R).
\end{equation*}
By Proposition~\ref{prop:global}, the right-hand side is bounded for $t\ge1$. Hence the local gradient energy is bounded on average, even though $u$ itself may grow linearly.
\end{rem}

The same identities give results for time-dependent profiles.

\begin{thm}\label{thm:profiles}
Let $N\ge1$, $q>1$ and $\mathtt a\in C([0,\infty))$. Then for every $R>0$ and $t\ge0$,
\begin{equation}\label{eq:profiles}
F_R(t)\ge F_R(0)+A(t)M(R)-C_q\Theta R^{N-q'}t .
\end{equation}
Consequently, the following hold.
\begin{enumerate}
\item[(i)] If $A(t)/t\to+\infty$ as $t\to\infty$, there are $\rho,c,t_1>0$ such that $\max_{|x|\le\rho}u(t,x)\ge c\,A(t)$ for $t\ge t_1$. If moreover $\mathtt a\ge0$, then $u(t,x)\le\sup u_0+\|\mathbf w_+\|_\infty A(t)$, so that $u$ grows exactly like $A(t)$ on a ball.
\item[(ii)] If $1<q\leq q_c$ and $\liminf_{t\to\infty}A(t)/t>0$, then $u$ grows linearly on a ball.
\end{enumerate}
\end{thm}

\begin{proof}
Dropping the gradient term in \eqref{eq:id2} and integrating, we obtain \eqref{eq:profiles}.

For (i), we fix $R$ with $M(R)\ge M_\infty/2$. Since $A(t)/t\to\infty$, we have $A(t)\to\infty$ and $t=o(A(t))$ as $t\to\infty$. By \eqref{eq:profiles}, we infer that $F_R(t)\ge\frac{M_\infty}{4}A(t)$ for large $t$. Applying \eqref{eq:avg} with $\rho=2R$, we get the lower bound with $c=M_\infty/(4\|\varphi_R\|_{L^1})$. The upper bound follows from the comparison argument of Proposition~\ref{prop:global}, using the supersolution $\sup u_0+\|\mathbf w_+\|_\infty A(t)$, which is valid when $\mathtt a\ge0$.

For (ii), let $\ell>0$ and $t_0>0$ be such that $A(t)\ge\ell t$ for $t\ge t_0$. We distinguish the subcritical and the critical cases.

\noindent\textit{Case 1: $1<q<q_c$.} Since $q'>N$, we may choose $R$ so large that $M(R)\ge M_\infty/2$ and $C_q\Theta R^{N-q'}\le\ell M_\infty/4$. Since $A(t)M(R)\ge\ell tM_\infty/2$ for $t\ge t_0$, we deduce from \eqref{eq:profiles} that $F_R(t)\ge F_R(0)+\ell M_\infty t/4$ for $t\ge t_0$. Using \eqref{eq:avg}, we conclude that
\begin{equation*}
\max_{|x|\le2R}u(t,x)\ge\inf u_0+\frac{\ell M_\infty}{4\|\varphi_R\|_{L^1}}\,t,\qquad t\ge t_0 .
\end{equation*}

\noindent\textit{Case 2: $q=q_c$, so that $N\ge2$.} We use Lemma~\ref{lem:identity}(b). Dropping the gradient term in \eqref{eq:id3}, integrating over $(\varepsilon,t)$ and letting $\varepsilon\to0$, we get
\begin{equation*}
G_R(t)\ge G_R(0)+A(t)\widetilde M(R)-\frac{C_q\Theta_N}{(\log R)^{N-1}}\,t .
\end{equation*}
Since $\widetilde M(R)\to M_\infty$ and $(\log R)^{1-N}\to0$, we can choose $R_0>1$ such that $\widetilde M(R_0)\ge M_\infty/2$ and $C_q\Theta_N(\log R_0)^{1-N}\le\ell M_\infty/4$. It follows that $G_{R_0}(t)\ge G_{R_0}(0)+\ell M_\infty t/4$ for $t\ge t_0$. Since $\operatorname{supp}\zeta_{R_0}\subset\overline{B_{R_0^2}}$, \eqref{eq:avg} yields
\begin{equation*}
\max_{|x|\le R_0^2}u(t,x)\ge\inf u_0+\frac{\ell M_\infty}{4\|\zeta_{R_0}\|_{L^1}}\,t,\qquad t\ge t_0 .
\end{equation*}

In both cases, the lower bound holds for $t\ge t_0$. On the bounded interval $[0,t_0]$, the solution is bounded from below by \eqref{eq:apriori}. Increasing the constant $C$ if necessary, we conclude that $u$ grows linearly on a ball in the sense of Definition~\ref{def:growup}.
\end{proof}

\begin{rem}
Part (i) holds in {every} dimension and for {every} $q>1$. It is the analogue, for gradient nonlinearities, of the fact that a forcing whose primitive grows superlinearly precludes global solutions of $u_t=\Delta u+|u|^p+\mathtt a(t)\mathbf w$ for every $p>1$ \cite{JKS}. For instance, $\mathtt a(t)=\log(e+t)$ leads to growth of order $t\log t$ for every $q>1$ and every $N$.
\end{rem}

For decaying profiles, the same identity with a time-dependent radius gives a partial result.

\begin{cor}\label{cor:decay}
Let $-1<\sigma<0$ and assume that $\liminf_{t\to\infty}t^{-1-\sigma}A(t)>0$. Suppose that $q'(1+\sigma)>N$, that is, $q<N/(N-1-\sigma)$. Set
\begin{equation*}
\theta:=\frac{-\sigma}{q'-N}>0,\qquad \kappa_\sigma:=1+\sigma-N\theta=\frac{(1+\sigma)q'-N}{q'-N}>0 .
\end{equation*}
Then there are $c,C,t_1>0$ such that $\max_{|x|\le Ct^\theta}u(t,x)\ge c\,t^{\kappa_\sigma}$ for $t\ge t_1$. In particular, $u$ grows up.
\end{cor}

\begin{proof}
Note first that $q'(1+\sigma)>N$ implies $q'>N$, so that $\theta$ is well defined. The equivalence with $q<N/(N-1-\sigma)$ follows from $N-1-\sigma>0$. Let $\ell>0$ and $t_0$ be such that $A(t)\ge\ell t^{1+\sigma}$ for $t\ge t_0$. We apply \eqref{eq:profiles} at time $t$ with $R=R(t):=\beta t^\theta$, where $\beta>0$ is chosen such that $C_q\Theta\beta^{N-q'}\le\ell M_\infty/4$. This is possible since $N-q'<0$. For $t$ large, we have $M(R(t))\ge M_\infty/2$. Moreover $C_q\Theta R(t)^{N-q'}t=C_q\Theta\beta^{N-q'}t^{1+\sigma}$, because $\theta(N-q')=\sigma$. Using $|F_{R}(0)|\le\|u_0\|_\infty|B_{2R}|$ and $|B_{2R(t)}|=|B_1|(2\beta)^Nt^{N\theta}$, we obtain
\begin{equation*}
F_{R(t)}(t)\ge-\|u_0\|_\infty|B_1|(2\beta)^Nt^{N\theta}+\frac{\ell M_\infty}{4}t^{1+\sigma}\ge\frac{\ell M_\infty}{8}t^{1+\sigma}
\end{equation*}
for large $t$, since $N\theta<1+\sigma$. In particular $F_{R(t)}(t)>0$, so that $\max_{|x|\le2R(t)}u(t,x)>0$ and $F_{R(t)}(t)\le|B_{2R(t)}|\max_{|x|\le2R(t)}u(t,x)$. Dividing by $|B_{2R(t)}|$, we get the claim with $C=2\beta$.
\end{proof}

\begin{rem}
Corollary~\ref{cor:decay} is not expected to be sharp. The scaling computation of Subsection~\ref{ss:heuristics}, applied to $t^\sigma\mathbf w$, suggests the threshold $(N-2\sigma)/(N-1-2\sigma)$. This number is strictly larger than $N/(N-1-\sigma)$ when $N\ge2$ and $\sigma<0$. Indeed, in terms of $q'$ the two conditions read $q'>N-2\sigma$ and $q'>N/(1+\sigma)$, and
\begin{equation*}
\frac{N}{1+\sigma}-(N-2\sigma)=\frac{\sigma(2+2\sigma-N)}{1+\sigma}>0 .
\end{equation*}
The loss presumably comes from discarding the gradient term in \eqref{eq:id2}, see Conjecture~\ref{conj:decay}.
\end{rem}

\section{The critical exponent}\label{sec:crit}

At the critical exponent $q=q_c$, the localisation error $C_q\Theta R^{N-q'}$ in \eqref{eq:id2} no longer vanishes as $R\to\infty$, since $q'=N$. To overcome this difficulty, we use the logarithmic cut-off functions of Lemma~\ref{logcutoff}, whose weighted gradient energy tends to zero as $R\to\infty$. This allows us to extend the linear growth of Theorem~\ref{thm:sub} to the critical case.

\begin{thm}\label{thm:crit}
Let $N\ge2$, $q=q_c=N/(N-1)$ and $\mathtt{a}\equiv1$. Then there exist constants $\rho,c>0$ and $C\ge0$ such that
\begin{equation}\label{eq:critconcl}
\max_{|x|\le\rho}u(t,x)\ge ct-C,\qquad t\ge0.
\end{equation}
In particular,
\begin{equation*}
\liminf_{t\to\infty}\,\frac{1}{t}\max_{|x|\le\rho}u(t,x)\ge c>0,
\end{equation*}
and every solution grows up.
\end{thm}

\begin{proof}
Let $G_R$ and $\widetilde M(R)$ be as in Lemma~\ref{lem:identity}(b). Since $\widetilde M(R)\to M_\infty>0$ and $C_q\Theta_N(\log R)^{1-N}\to0$ as $R\to\infty$, we can choose $R_0>1$ such that
\begin{equation*}
\widetilde M(R_0)-\frac{C_q\Theta_N}{(\log R_0)^{N-1}}\ge\frac{M_\infty}{2}.
\end{equation*}
Dropping the nonnegative gradient term in \eqref{eq:id3} with $\mathtt a\equiv1$, we obtain $G_{R_0}'(t)\ge M_\infty/2$ for $t>0$. Integrating over $(\varepsilon,t)$ and letting $\varepsilon\to0$, which is allowed since $G_{R_0}$ is continuous at $t=0$, we get
\begin{equation}\label{EQ1}
G_{R_0}(t)\ge G_{R_0}(0)+\frac{M_\infty}{2}\,t,\qquad t\ge0.
\end{equation}
Set $J_{R_0}:=\int_{\R^N}\zeta_{R_0}\,\dd x$. Since $0\le\zeta_{R_0}\le1$, $\zeta_{R_0}=1$ on $B_{R_0}$ and $\operatorname{supp}\zeta_{R_0}\subset\overline{B_{R_0^2}}$, we have $0<J_{R_0}<\infty$. By \eqref{eq:avg}, we get $G_{R_0}(t)\le J_{R_0}\max_{|x|\le R_0^2}u(t,x)$ and $G_{R_0}(0)\ge-\|u_0\|_\infty J_{R_0}$. Combining these estimates with \eqref{EQ1} and dividing by $J_{R_0}$, we conclude that
\begin{equation*}
\max_{|x|\le R_0^2}u(t,x)\ge\frac{M_\infty}{2J_{R_0}}\,t-\|u_0\|_\infty,\qquad t\ge0 .
\end{equation*}
This is \eqref{eq:critconcl} with $\rho=R_0^2$, $c=M_\infty/(2J_{R_0})$ and $C=\|u_0\|_\infty$. Dividing \eqref{eq:critconcl} by $t$ and letting $t\to\infty$, we obtain the lower bound for the $\liminf$. This completes the proof.
\end{proof}

\begin{rem}
Theorem~\ref{thm:crit} is the special case $\mathtt a\equiv1$ of Theorem~\ref{thm:profiles}(ii). We state it separately because it settles the critical exponent for stationary forcing. The linear lower bound holds on a fixed ball and for all $t\ge0$. No pointwise sign condition on $\mathbf w$ is required, apart from the positivity of its total mass. When $N=2$, the critical exponent is $q_c=2$, and Theorem~\ref{thm:quadratic} below provides an independent spectral lower bound for the growth rate.
\end{rem}

\section{The plane: grow-up for every exponent}\label{sec:plane}

In dimension two, $q_c=2$. For $q>2$, the gradient argument of Sections~\ref{sec:sub} and \ref{sec:crit} fails. Indeed, $q>2$ means $q'<2$, and the $q'$-capacity of a ball relative to infinity is positive in $\R^2$ when $q'<2$. Hence no cut-off function equal to $1$ on a large ball can make the weighted energy in \eqref{eq:id2} small. We use instead identity \eqref{eq:id1} with the logarithmic cut-off functions of Lemma~\ref{lem:log}, together with the $L^\infty$ bound of Proposition~\ref{prop:global}. We then prove a sharper, local statement by comparison with the linear heat equation.

\begin{thm}\label{thm:plane}
Let $N=2$, $q>1$ and $\mathtt a\equiv1$. Then every solution grows up. More precisely, $K(t):=\sup_{0\le s\le t}\|u(s)\|_\infty$ satisfies
\begin{equation*}
K(t)\ge c\log t\qquad\text{for all large }t,
\end{equation*}
where $c>0$ depends only on $M_\infty$ and on the function $\eta$ of Lemma~\ref{lem:log}.
\end{thm}

\begin{proof}
We repeat the computation of Lemma~\ref{lem:identity}(a) with $\ell_R$ in place of $\varphi_R$. For $R\ge2$, we set $\widehat G_R(t):=\int u(t)\ell_R$ and $\widehat M(R):=\int\mathbf w\ell_R$. By dominated convergence, $\widehat M(R)\to M_\infty$ as $R\to\infty$. Identity \eqref{eq:id1} reads $\widehat G_R'=\int u\Delta\ell_R+\int|\nabla u|^q\ell_R+\widehat M(R)$. Dropping the gradient term and using Lemma~\ref{lem:log}, we get
\begin{equation*}
\widehat G_R'(s)\ge\widehat M(R)-\|u(s)\|_\infty\frac{2\pi\|\eta''\|_{L^1}}{\log R}.
\end{equation*}
Integrating over $(0,t)$, we obtain $\widehat G_R(t)-\widehat G_R(0)\ge t\widehat M(R)-tK(t)\,2\pi\|\eta''\|_{L^1}/\log R$. On the other hand, $|\widehat G_R(t)-\widehat G_R(0)|\le2K(t)\,\pi R^4$, because $0\le\ell_R\le1$ and $\operatorname{supp}\ell_R\subset\overline{B_{R^2}}$. Hence, for $R$ so large that $\widehat M(R)\ge M_\infty/2$, we infer that
\begin{equation*}
K(t)\Big(2\pi R^4+\frac{2\pi\|\eta''\|_{L^1}\,t}{\log R}\Big)\ge\frac{M_\infty}{2}\,t .
\end{equation*}
We choose $R=R(t):=(t/\log t)^{1/4}$, which tends to infinity with $t$. Then $2\pi R^4=2\pi t/\log t$ and $\log R=\frac14(\log t-\log\log t)\ge\frac18\log t$ for large $t$. We conclude that
\begin{equation*}
K(t)\ge\frac{M_\infty}{2}\,\frac{\log t}{2\pi+16\pi\|\eta''\|_{L^1}}\qquad\text{for large }t. \qedhere
\end{equation*}
\end{proof}

\begin{rem}
Theorem~\ref{thm:plane} uses neither the size of $q$ nor the gradient term. It reflects only the recurrence of planar Brownian motion, which is encoded in Lemma~\ref{lem:log}. A positive mass supplied at a constant rate cannot be dispersed to infinity in two dimensions. For $q<2$, Theorem~\ref{thm:sub} gives linear growth, and Theorems~\ref{thm:crit} and \ref{thm:quadratic} give linear growth for $q=2$. We conjecture that growth is linear for every $q>1$ when $N=2$ (Conjecture~\ref{conj:ergodic}). For $N=1$, linear growth for every $q>1$ is already contained in Theorem~\ref{thm:sub}, since $q_c=\infty$.
\end{rem}

The next result shows that the divergence in the plane is local, uniform and has an explicit logarithmic rate.

\begin{thm}\label{thm:plane-local}
Let $N=2$, $q>1$ and $\mathtt a\equiv1$, and let $u_0$ and $\mathbf w$ satisfy the assumptions of Proposition~\ref{prop:global} and \eqref{eq:mass}. Then, for every $R>0$,
\begin{equation*}
\liminf_{t\to\infty}\frac{1}{\log t}\inf_{|x|\le R}u(t,x)\geq\frac{M_\infty}{4\pi}>0.
\end{equation*}
In particular, $u(t,x)\to+\infty$ locally uniformly in $\R^2$ as $t\to\infty$, and Theorem~\ref{thm:plane} holds with any constant $c<M_\infty/(4\pi)$.
\end{thm}

\begin{proof}
Let $U$ be the solution of the linear problem $U_t=\Delta U+\mathbf w$ with $U(0)=u_0$, that is,
\begin{equation*}
U(t)=e^{t\Delta}u_0+\int_0^te^{s\Delta}\mathbf w\,\dd s .
\end{equation*}
Since $u_0,\mathbf w\in W^{1,\infty}(\R^2)$ and $\|D^2e^{s\Delta}\mathbf w\|_\infty\le Cs^{-1/2}\|\nabla\mathbf w\|_\infty$ is integrable near $s=0$, the function $U$ is a classical solution, and $U$ and $\nabla U$ are bounded on $(0,T]\times\R^2$ for every $T>0$, with $|U(t)|\le\|u_0\|_\infty+t\|\mathbf w\|_\infty$. The function $z=U-u$ satisfies $z_t-\Delta z=-|\nabla u|^q\le0$, $z(0)=0$, and it is bounded from above on $[0,T]\times\R^2$. Lemma~\ref{lem:maxprinciple} with $b\equiv0$ gives $u\ge U$ on $[0,\infty)\times\R^2$.

We now estimate $U$ from below on $\overline{B_R}$. Clearly $e^{t\Delta}u_0\ge\inf u_0$. For $s>0$ and $|x|\le R$, we have
\begin{equation*}
4\pi s\,(e^{s\Delta}\mathbf w)(x)-M_\infty=\int_{\R^2}\big(e^{-|x-y|^2/(4s)}-1\big)\mathbf w(y)\,\dd y .
\end{equation*}
Since $0\le1-e^{-a}\le\min\{1,a\}$ for $a\ge0$ and $|x-y|\le R+|y|$, we get
\begin{equation*}
\big|4\pi s\,(e^{s\Delta}\mathbf w)(x)-M_\infty\big|\le\varepsilon_R(s):=\int_{\R^2}\min\Big\{1,\frac{(R+|y|)^2}{4s}\Big\}|\mathbf w(y)|\,\dd y .
\end{equation*}
By dominated convergence, $\varepsilon_R(s)\to0$ as $s\to\infty$. Fix $\varepsilon\in(0,M_\infty)$ and choose $s_0\ge1$ such that $\varepsilon_R(s)\le\varepsilon$ for $s\ge s_0$. Using $|e^{s\Delta}\mathbf w|\le\|\mathbf w\|_\infty$ for $s\le s_0$, we deduce that, for $t\ge s_0$ and $|x|\le R$,
\begin{equation*}
\int_0^t(e^{s\Delta}\mathbf w)(x)\,\dd s\ge-s_0\|\mathbf w\|_\infty+\frac{M_\infty-\varepsilon}{4\pi}\int_{s_0}^t\frac{\dd s}{s}=-s_0\|\mathbf w\|_\infty+\frac{M_\infty-\varepsilon}{4\pi}\log\frac{t}{s_0}.
\end{equation*}
It follows that
\begin{equation*}
\inf_{|x|\le R}u(t,x)\ge\inf u_0-s_0\|\mathbf w\|_\infty+\frac{M_\infty-\varepsilon}{4\pi}\log\frac{t}{s_0},\qquad t\ge s_0 .
\end{equation*}
Dividing by $\log t$ and letting $t\to\infty$, we obtain that the $\liminf$ is at least $(M_\infty-\varepsilon)/(4\pi)$. Since $\varepsilon$ is arbitrary, the claim follows. The last assertion follows from $K(t)\ge u(t,0)$.
\end{proof}

\begin{rem}
The proof of Theorem~\ref{thm:plane-local} only uses the inequality $u\ge U$, which holds in every dimension and for every $q>1$. When $N=1$, the same computation gives $\inf_{|x|\le R}u(t,x)\ge c\sqrt t$ for large $t$, which is weaker than Theorem~\ref{thm:sub}. When $N\ge3$, it only gives a lower bound that is uniform in time, because $\int_1^\infty s^{-N/2}\,\dd s<\infty$. The constant $M_\infty/(4\pi)$ is the exact rate for the linear problem, which corresponds to dropping the gradient term. Any improvement beyond the logarithmic rate for $q>2$ must therefore use the nonlinearity.
\end{rem}

\section{The quadratic case}\label{sec:quadratic}

For $q=2$, the Cole--Hopf transformation \eqref{eq:colehopf} gives a complete and quantitative picture.

\begin{thm}\label{thm:quadratic}
Let $q=2$, $\mathtt a\equiv1$ and $\mathbf w\in W^{1,\infty}(\R^N)$. Here we do not need $\mathbf w\in L^1$. Then, for every $R>0$ and $t\ge0$,
\begin{equation}\label{eq:quadlower}
\max_{|x|\le R}u(t,x)\ge\inf u_0-\lambda_R\,t .
\end{equation}
Consequently, if $\lambda_0<0$, where $\lambda_0$ is defined by \eqref{eq:lambda0}, then
\begin{equation*}
\liminf_{t\to\infty}\frac1t\sup_{x\in\R^N}u(t,x)\ge-\lambda_0>0 .
\end{equation*}
\end{thm}

\begin{proof}
Let $v=e^u$. Then $v$ solves \eqref{eq:colehopf} with $\mathtt a\equiv1$, and $v(0)\ge e^{\inf u_0}>0$. Let $H(t):=\int_{B_R}v(t)\psi_R$, where $\psi_R$ is the normalised principal eigenfunction introduced in Subsection~\ref{ss:colehopf}. By Green's formula, using $\psi_R=0$ and $\partial_\nu\psi_R\le0$ on $\partial B_R$ together with $v>0$, we obtain
\begin{equation*}
\begin{aligned}
H'(t)&=\int_{B_R}(\Delta v+\mathbf wv)\psi_R\\
&=\int_{B_R}v\,(\Delta\psi_R+\mathbf w\psi_R)-\int_{\partial B_R}v\,\partial_\nu\psi_R\\
&\ge-\lambda_R H(t).
\end{aligned}
\end{equation*}
Hence $H(t)\ge e^{-\lambda_Rt}H(0)\ge e^{-\lambda_Rt}e^{\inf u_0}$. Since $\int_{B_R}\psi_R=1$ and $\psi_R\ge0$, we also have $H(t)\le\exp\big(\max_{|x|\le R}u(t,x)\big)$. Taking logarithms, we get \eqref{eq:quadlower}. If $\lambda_0<0$, then for every $\varepsilon\in(0,-\lambda_0)$ there is $R$ with $\lambda_R\le\lambda_0+\varepsilon<0$. By \eqref{eq:quadlower}, the $\liminf$ is at least $-\lambda_0-\varepsilon$, and the second claim follows.
\end{proof}

\begin{cor}\label{cor:quadratic}
Let $q=2$ and $\mathtt a\equiv1$.
\begin{enumerate}
\item[(a)] If $N\in\{1,2\}$, $\mathbf w\in L^1\cap W^{1,\infty}$, $\int\mathbf w\ge0$ and $\mathbf w\not\equiv0$, then $\lambda_0<0$. Hence every solution grows at least linearly, with rate at least $-\lambda_0$. In particular, this holds under \eqref{eq:mass}, and also for nontrivial forcing with zero mass.
\item[(b)] If $\mathbf w_1\in W^{1,\infty}$, $\mathbf w_1\ge0$ and $\mathbf w_1\not\equiv0$, then $\lambda_0<0$ for $\mathbf w=\kappa\mathbf w_1$ with $\kappa$ large enough. Hence solutions grow linearly for large forcing. This is relevant when $N\ge3$.
\item[(c)] If $N\ge3$ and $\mathbf w(x)\le\frac{(N-2)^2}{4}|x|^{-2}$, then $\lambda_0\ge0$.
\end{enumerate}
\end{cor}

\begin{proof}
(a) Let $\eta$ be as in Lemma~\ref{lem:log}. Let $f_R(x)=\eta(|x|/R)$ if $N=1$, and $f_R=\ell_R$ if $N=2$. In both cases, $f_R=1$ on $B_R$ and $\int|\nabla f_R|^2\to0$ as $R\to\infty$. Indeed, for $N=1$ a direct computation gives $\int|f_R'|^2=2R^{-1}\|\eta'\|_{L^2}^2$, and for $N=2$ we use Lemma~\ref{lem:log}. By dominated convergence, $\int\mathbf wf_R^2\to\int\mathbf w$. If $\int\mathbf w>0$, we deduce that $\int(|\nabla f_R|^2-\mathbf wf_R^2)<0$ for $R$ large, and $\lambda_0<0$ by \eqref{eq:lambda0}.

Assume now that $\int\mathbf w=0$ and $\mathbf w\not\equiv0$. Since $\mathbf w$ is continuous and not identically zero, there is $g\in C_c^\infty(\R^N)$ with $\int\mathbf wg>0$. Let $R$ be so large that $\operatorname{supp}g\subset B_R$, where $f_R=1$ and $\nabla f_R=0$. For $f=f_R(1+\varepsilon g)=f_R+\varepsilon g$, we have $\int\nabla f_R\cdot\nabla g=0$ and $f^2=f_R^2+2\varepsilon g+\varepsilon^2g^2$. This gives
\begin{equation*}
\int\big(|\nabla f|^2-\mathbf wf^2\big)=\int|\nabla f_R|^2-\int\mathbf wf_R^2-2\varepsilon\int\mathbf wg+\varepsilon^2\int\big(|\nabla g|^2-\mathbf wg^2\big).
\end{equation*}
As $R\to\infty$, the first two terms tend to $0$ and to $-\int\mathbf w=0$, respectively. Choosing first $\varepsilon$ small and then $R$ large, we make the right-hand side negative, and we conclude that $\lambda_0<0$.

(b) We fix $f\in C_c^\infty(\R^N)$ with $\int\mathbf w_1f^2>0$. Then $\int|\nabla f|^2-\kappa\int\mathbf w_1f^2<0$ for $\kappa$ large.

(c) By Hardy's inequality $\int|\nabla f|^2\ge\frac{(N-2)^2}{4}\int|x|^{-2}f^2$, valid for $f\in C_c^\infty(\R^N)$ and $N\ge3$ \cite{HLP52, BEL15}, the quadratic form in \eqref{eq:lambda0} is nonnegative.
\end{proof}

\begin{rem}[The exact growth rate]\label{rem:exact}
The upper bound $\sup u(t)\le\sup u_0+t\|\mathbf w_+\|_\infty$ of Proposition~\ref{prop:global} is not sharp. A sharp upper bound follows from known results on Schr\"odinger semigroups, which we now recall. For bounded $\mathbf w$, the spectrum of $-\Delta-\mathbf w$ in $L^p(\R^N)$ does not depend on $p\in[1,\infty]$ \cite{HV86, Sim82}. For positive $C_0$-semigroups on $L^1$, the growth bound coincides with the spectral bound \cite{Wei95}. Since the kernel of $e^{t(\Delta+\mathbf w)}$ is symmetric, we have $\|e^{t(\Delta+\mathbf w)}\|_{L^\infty\to L^\infty}=\|e^{t(\Delta+\mathbf w)}\|_{L^1\to L^1}$. We infer that $t^{-1}\log\|e^{t(\Delta+\mathbf w)}\|_{L^\infty\to L^\infty}\to-\lambda_0$ as $t\to\infty$, see also \cite{Sim81}. The classical solution $v$ of \eqref{eq:colehopf} coincides with $e^{t(\Delta+\mathbf w)}e^{u_0}$, by uniqueness of bounded classical solutions on finite time intervals. Since $e^{u_0}\le e^{\sup u_0}$, it follows that $\limsup_{t\to\infty}t^{-1}\sup u(t)\le-\lambda_0$. Combined with Theorem~\ref{thm:quadratic}, this gives
\begin{equation*}
\lim_{t\to\infty}\frac1t\sup_{x\in\R^N}u(t,x)=-\lambda_0\qquad\text{whenever }\lambda_0<0 .
\end{equation*}
We state this identity as a consequence of the cited results, which we do not reprove here. When $\lambda_0\ge0$, the solution may remain bounded, see Theorem~\ref{thm:bounded} for small $\mathbf w$. The borderline case $\lambda_0=0$ is related to the criticality theory of Schr\"odinger operators \cite{Mur86, Pin07}.
\end{rem}

\begin{rem}
Corollary~\ref{cor:quadratic}(a) is the weak coupling phenomenon. In dimensions one and two, every potential with positive mass creates a bound state \cite{Sim76}. Under \eqref{eq:mass}, $\mathbf w$ is Lipschitz and integrable, hence $\mathbf w(x)\to0$ as $|x|\to\infty$, and the essential spectrum of $-\Delta-\mathbf w$ is $[0,\infty)$ \cite{RS78}. Consequently, $\lambda_0<0$ is an eigenvalue. The dichotomy is consistent with the exponent $q_c=N/(N-1)$. For $N=2$, we have $q=2=q_c$, and linear growth holds at the critical exponent. For $N\ge3$, we have $q=2>q_c$, and the outcome depends on the size of $\mathbf w$, as shown by parts (b) and (c) and by Theorem~\ref{thm:bounded}.
\end{rem}

\section{Bounded solutions in the supercritical range and a large-forcing criterion}\label{sec:bounded}

We first construct explicit stationary supersolutions for $q>q_c$ and $N\ge3$. The exponent $q_c$ appears exactly as the condition for such supersolutions to exist.

\begin{lem}\label{lem:supersolution}
Let $N\ge3$ and $q>q_c=N/(N-1)$. Then the interval $\big(\max\{0,\frac{2-q}{q-1}\},\,N-2\big)$ is nonempty. Fix $\gamma$ in this interval, and let
\begin{equation*}
0<\delta\le\delta_*:=\frac1\gamma\Big(\frac{N-2-\gamma}{2}\Big)^{\frac1{q-1}},\qquad \overline U(x):=\delta\jb{x}^{-\gamma}.
\end{equation*}
Then
\begin{equation*}
-\Delta\overline U-|\nabla\overline U|^q\ge c_\delta\jb{x}^{-\gamma-2}\quad\text{in }\R^N,\qquad c_\delta:=\frac{\delta\gamma(N-2-\gamma)}{2}.
\end{equation*}
\end{lem}

\begin{proof}
The inequality $\frac{2-q}{q-1}<N-2$ is equivalent to $2-q<(N-2)(q-1)$, that is, to $q(N-1)>N$. Hence the interval is nonempty precisely because $q>q_c$. A direct computation gives
\begin{equation*}
-\Delta\jb{x}^{-\gamma}=\gamma\jb{x}^{-\gamma-4}\big(N+(N-2-\gamma)|x|^2\big)\ge\gamma(N-2-\gamma)\jb{x}^{-\gamma-2},
\end{equation*}
because $N\ge N-2-\gamma>0$. We also have $|\nabla\overline U|=\delta\gamma|x|\jb{x}^{-\gamma-2}\le\delta\gamma\jb{x}^{-\gamma-1}$. Since $\gamma>\frac{2-q}{q-1}$ means $q(\gamma+1)>\gamma+2$, we get $|\nabla\overline U|^q\le(\delta\gamma)^q\jb{x}^{-\gamma-2}$. Finally, $(\delta\gamma)^q=\delta\gamma(\delta\gamma)^{q-1}\le\delta\gamma(N-2-\gamma)/2$ by the choice of $\delta$. Combining these estimates, we obtain the claim.
\end{proof}

\begin{thm}\label{thm:bounded}
Let $N\ge3$, $q>q_c$, and let $\gamma$, $\delta$ and $c_\delta$ be as in Lemma~\ref{lem:supersolution}. Assume that $\mathtt a\in C([0,\infty))$ satisfies $0\le\mathtt a\le1$, that $\mathbf w\in W^{1,\infty}(\R^N)$ satisfies $0\le\mathbf w\le c_\delta\jb{x}^{-\gamma-2}$, and that $u_0\in W^{1,\infty}(\R^N)$ satisfies $0\le u_0\le\delta\jb{x}^{-\gamma}$. Then the solution of \eqref{eq:main} satisfies
\begin{equation*}
0\le u(t,x)\le\delta\jb{x}^{-\gamma}\qquad\text{for all }(t,x)\in[0,\infty)\times\R^N .
\end{equation*}
The admissible forcing profiles include functions satisfying \eqref{eq:mass}, for instance $\mathbf w=c_\delta\jb{x}^{-N-1}$.
\end{thm}

\begin{proof}
By Lemma~\ref{lem:supersolution}, $\overline U$ satisfies $\partial_t\overline U-\Delta\overline U-|\nabla\overline U|^q-\mathtt a\mathbf w\ge c_\delta\jb{x}^{-\gamma-2}-\mathbf w\ge0$, since $0\le\mathtt a\le1$ and $\mathbf w\ge0$. The function $0$ is a subsolution, since $-\mathtt a\mathbf w\le0$. Both functions are bounded with bounded gradients, and $0\le u_0\le\overline U$. On every interval $[0,T]$, Lemma~\ref{lem:comparison} with $g=\mathtt a\mathbf w$ gives $0\le u\le\overline U$. Finally, $\jb{x}^{-N-1}\le\jb{x}^{-\gamma-2}$ because $\gamma+2<N<N+1$, and $\jb{x}^{-N-1}$ is integrable and Lipschitz.
\end{proof}

\begin{prop}\label{prop:monotone}
Under the assumptions of Theorem~\ref{thm:bounded} with $\mathtt a\equiv1$ and $u_0\equiv0$, the map $t\mapsto u(t,x)$ is nondecreasing for every $x$. Consequently, $u(t,\cdot)$ converges pointwise, as $t\to\infty$, to a function $U$ with $0\le U\le\delta\jb{x}^{-\gamma}$.
\end{prop}

\begin{proof}
Fix $h>0$. Since the equation is autonomous when $\mathtt a\equiv1$, the function $u_h(t,x):=u(t+h,x)$ solves \eqref{eq:main} with initial value $u(h)$, and $u(h)\ge0=u_0$ by Theorem~\ref{thm:bounded}. By Lemma~\ref{lem:comparison}, we get $u_h\ge u$. The monotonicity and the bound of Theorem~\ref{thm:bounded} give the convergence.
\end{proof}

\begin{rem}
We expect $U$ to be the minimal nonnegative classical solution of $-\Delta U=|\nabla U|^q+\mathbf w$. This would follow from interior gradient estimates for bounded solutions that are uniform in time. Such estimates are available for the stationary problem \cite{Lio85}. We do not pursue this here, and we do not claim it.
\end{rem}

Theorem~\ref{thm:bounded} requires small forcing. Conversely, a scale-invariant largeness condition forces linear growth in every dimension and for every $q>1$.

\begin{thm}\label{thm:large}
Let $N\ge1$, $q>1$ and $\mathtt a\equiv1$, and let $C_q=\frac{q-1}{q}\big(\frac2q\big)^{\frac1{q-1}}$ be the constant in \eqref{eq:id2}. Assume that
\begin{equation}\label{eq:large}
R^{q'-N}\int_{\R^N}\mathbf w\,\varphi_R\,\dd x>C_q\Theta\qquad\text{for some }R>0.
\end{equation}
Then $F_R(t)\ge F_R(0)+\varepsilon t$ for all $t\ge0$, with $\varepsilon:=M(R)-C_q\Theta R^{N-q'}>0$. In particular,
\begin{equation*}
\max_{|x|\le2R}u(t,x)\ge\inf u_0+\frac{\varepsilon}{\|\varphi_R\|_{L^1}}\,t,\qquad t\ge0,
\end{equation*}
so that $u$ grows linearly on $B_{2R}$.
\end{thm}

\begin{proof}
Condition \eqref{eq:large} means exactly that $\varepsilon>0$. Dropping the gradient term in \eqref{eq:id2}, we obtain $F_R'\ge\varepsilon$ on $(0,\infty)$, and the lower bound for $F_R$ follows by integration. The lower bound for the maximum follows from \eqref{eq:avg}.
\end{proof}

\begin{rem}
Replacing $\mathbf w$ by $\mathbf w_\lambda:=\lambda^{q'}\mathbf w(\lambda\,\cdot)$, which corresponds to the scaling of Subsection~\ref{ss:heuristics}, we get $R^{q'-N}\int\mathbf w_\lambda\varphi_R=(\lambda R)^{q'-N}M(\lambda R)$. Hence the family of conditions \eqref{eq:large} is invariant under this scaling. For $q<q_c$, the quantity $R^{q'-N}M(R)$ tends to $+\infty$ as $R\to\infty$ whenever $\int\mathbf w>0$, which gives another proof of the linear growth in Theorem~\ref{thm:sub}. For $q>q_c$ and $\mathbf w\in L^1$, it tends to $0$ as $R\to\infty$, and also as $R\to0$ since $|M(R)|\le\|\mathbf w\|_\infty|B_{2R}|$. Condition \eqref{eq:large} is then a genuine largeness condition at some finite scale. It is satisfied, for instance, by $\kappa\mathbf w_1$ with $\kappa$ large, as soon as $\int\mathbf w_1\varphi_R>0$ for some $R$. Theorems~\ref{thm:bounded} and \ref{thm:large} together show that, for $q>q_c$ and $N\ge3$, the long-time behaviour depends on the size of the forcing, just as the Fujita dichotomy depends on the size of the data above the Fujita exponent.
\end{rem}

\section{Summary, comments and open problems}\label{sec:open}

\subsection{Summary}
Table~\ref{tab:summary} summarises what is proved when $\int\mathbf w>0$.

\begin{table}[ht]
\centering
\small
\renewcommand{\arraystretch}{1.3}
\begin{tabular}{|l|l|l|}
\hline
Range & Proved behaviour & Reference\\
\hline
$N=1$, $q>1$ & linear growth on a ball, all data & Thm.~\ref{thm:sub}\\
$N\ge2$, $1<q<q_c$ & linear growth on a ball, all data & Thm.~\ref{thm:sub}\\
$N\ge2$, $q=q_c$ & linear growth on a ball, all data & Thm.~\ref{thm:crit}\\
$N=2$, $q=2$ & linear growth, exact rate $-\lambda_0>0$ & Cor.~\ref{cor:quadratic}, Rem.~\ref{rem:exact}\\
$N=2$, $q>2$ & $\inf_{|x|\le R}u(t)\ge\big(\frac{M_\infty}{4\pi}-o(1)\big)\log t$ & Thms.~\ref{thm:plane}, \ref{thm:plane-local}\\
$N\ge3$, $q>q_c$, small data & bounded, $0\le u\le\delta\jb{x}^{-\gamma}$ & Thm.~\ref{thm:bounded}\\
any $N$, $q>1$, large forcing \eqref{eq:large} & linear growth on a ball & Thm.~\ref{thm:large}\\
any $N$, $q>1$, $A(t)/t\to\infty$ & growth of order $A(t)$ & Thm.~\ref{thm:profiles}(i)\\
any $N$, $1<q\le q_c$, $\liminf_{t\to\infty}A(t)/t>0$ & linear growth on a ball & Thm.~\ref{thm:profiles}(ii)\\
$A(t)\gtrsim t^{1+\sigma}$, $\sigma\in(-1,0)$, $q'(1+\sigma)>N$ & power growth $t^{\kappa_\sigma}$ & Cor.~\ref{cor:decay}\\
\hline
\end{tabular}
\medskip
\caption{Proved results. In the first seven rows $\mathtt a\equiv1$. All solutions are global (Proposition~\ref{prop:global}), and $q_c=N/(N-1)$ for $N\ge2$.}
\label{tab:summary}
\end{table}

\subsection{Comments}
The exponent $q_c=N/(N-1)$ plays for $u_t=\Delta u+|\nabla u|^q+\mathbf w$ the role that the Bandle--Levine--Zhang exponent $N/(N-2)$ plays for $u_t=\Delta u+|u|^p+\mathbf w$. The mechanism is different, however. Finite-time blow-up is replaced by grow-up, and the dimension enters through the capacity of balls relative to infinity, that is, through the parabolicity of $\R^N$ with respect to the $q'$-energy \cite{Gri99}. The distinction between the subcritical and the critical cases is reflected in the energy cost of localisation. When $q'>N$, the standard cut-off functions yield an error of order $R^{N-q'}$, which vanishes as $R\to\infty$. At the endpoint $q'=N$, the logarithmic cut-off functions of Lemma~\ref{logcutoff} yield an error of order $(\log R)^{1-N}$, which also tends to zero. Thus, in both cases, the positive total mass of the forcing can be detected with an arbitrarily small localisation error. In the quadratic case, this is the classical fact that $-\Delta-\mathbf w$ has a bound state in dimensions one and two \cite{Sim76}. The plane is special, since Theorems~\ref{thm:plane} and \ref{thm:plane-local} hold for every $q$. The reason is that the relevant capacity there is the $2$-capacity of the Laplacian, and not the $q'$-capacity attached to the gradient term.

\subsection{Conjectures and open problems}
The following statements are {not} proved. They are supported by the special cases and heuristics indicated.

\begin{conj}[Ergodic growth rate]\label{conj:ergodic}
Assume either $N\le2$, or $N\ge3$ and $q\le q_c$. Then there is $\lambda^*=\lambda^*(q,\mathbf w)>0$ such that $u(t,x)/t\to\lambda^*$ locally uniformly as $t\to\infty$, for every initial datum $u_0\in W^{1,\infty}$. The constant $\lambda^*$ is the generalised principal eigenvalue (ergodic constant) of the problem $-\Delta\phi-|\nabla\phi|^q-\mathbf w=-\lambda$ in $\R^N$.
\end{conj}
For $q=2$, the growth rate of $\sup u$ is $-\lambda_0$ (Remark~\ref{rem:exact}). The local uniform convergence of $u/t$ is then expected from the exponential decay of the ground state of $-\Delta-\mathbf w$. For coercive potentials, results of this type are known \cite{BS01, BPT10, Tch10, KS06, Ich12, Ich15, BQR21, QR24}. The non-coercive, decaying case relevant here is related to the criticality theory of \cite{Ich13}, and to the analysis of the generalised principal eigenvalue for superquadratic Hamiltonians and potentials vanishing at infinity in \cite{CI16}. In particular, the conjecture would give linear growth for $N=2$ and $q>2$, improving the logarithmic bounds of Theorems~\ref{thm:plane} and \ref{thm:plane-local}.

\begin{prob}[Threshold in the supercritical range]\label{prob:threshold}
Let $N\ge3$ and $q>q_c$. Characterise the forcing terms $\mathbf w$ for which bounded solutions exist. For $q=2$, Theorem~\ref{thm:quadratic} shows that $\lambda_0\ge0$ is necessary, and the answer is expected to be $\lambda_0\ge0$, with delicate behaviour in the critical case $\lambda_0=0$ \cite{Mur86, Pin07}. Theorems~\ref{thm:bounded} and \ref{thm:large} give sufficient conditions on either side. Is there a single scale-invariant functional of $\mathbf w$, in the spirit of \eqref{eq:large}, that decides the alternative?
\end{prob}

\begin{conj}[Decaying forcing]\label{conj:decay}
Let $\mathtt a(t)=(1+t)^\sigma$ with $-1<\sigma<0$, and let $N\ge2$. Solutions grow up for $q<q^*(\sigma):=(N-2\sigma)/(N-1-2\sigma)$. For $q>q^*(\sigma)$, bounded solutions exist for small data.
\end{conj}
The exponent $q^*(\sigma)$ is the $L^1$-critical exponent for the forcing under the scaling of Subsection~\ref{ss:heuristics}. Corollary~\ref{cor:decay} proves grow-up only for $q<N/(N-1-\sigma)$, which is smaller than $q^*(\sigma)$.

\begin{prob}[Zero and negative mass]\label{prob:zeromass}
For $q=2$ and $N\le2$, Corollary~\ref{cor:quadratic}(a) shows that nontrivial forcing with zero mass still produces linear growth. Does the same hold for $q\ne2$? What happens in higher dimensions, where forcing with zero mass is expected to be weaker than any forcing with positive mass?
\end{prob}

\begin{prob}[Other diffusions and weights]\label{prob:other}
For a diffusion operator of order $\beta\in(1,2]$, for instance $-(-\Delta)^{\beta/2}$, the same scaling computation suggests the critical exponent $q_c(\beta)=N/(N+1-\beta)$. Prove the corresponding growth and boundedness results. Treat also weighted gradient terms $|x|^\alpha|\nabla u|^q$, combined nonlinearities $|u|^p+|\nabla u|^q$, gradient absorption $-|\nabla u|^q$ with positive forcing, and inverse-square potentials, for which Hardy-type inequalities should replace Corollary~\ref{cor:quadratic}(c).
\end{prob}

\begin{prob}[Rough data]
Extend the results to initial data and forcing terms in Lebesgue spaces, using the local theory of \cite{BASW02, BL99} and the notion of viscosity solutions \cite{CIL92}. In particular, determine the minimal regularity of $\mathbf w$ under which grow-up persists.
\end{prob}

\section*{Declarations}
\noindent\textbf{Funding.} This research received no external funding.

\noindent\textbf{Conflict of interest.} The authors declare that they have no conflicts of interest.

\noindent\textbf{Data availability.} No datasets were generated or analyzed during the current study.

\end{document}